\documentclass[12pt,draft]{article}

\usepackage{amsfonts}
\usepackage{amsmath}
\usepackage{amssymb}
\usepackage{amscd}
\usepackage{amsthm}
\usepackage{indentfirst}
\usepackage[vmargin=4cm]{geometry}

\newtheorem{thm}{Theorem}[section]
\newtheorem{lemma}[thm]{Lemma}
\newtheorem{prop}[thm]{Proposition}
\newtheorem{coroll}[thm]{Corollary}

\theoremstyle{definition}

\newtheorem{defin}[thm]{Definition}
\newtheorem{rem}[thm]{Remark}

\newtheorem*{acknow}{Acknowledgements}
\newtheorem*{prf}{Proof}

\newcommand{\R}{{\mathbb{R}}}

\newcommand{\T}{{\mathbb{T}}}
\newcommand{\Z}{{\mathbb{Z}}}

\newcommand{\cA}{{\mathcal{A}}}

\newcommand{\fc}{{:\ }}

\newcommand{\ol}{\overline}

\newcommand{\tb}{\textbf}

\DeclareMathOperator{\im}{im}

\DeclareMathOperator{\const}{const}

\DeclareMathOperator{\Int}{Int}

\DeclareMathOperator{\ed}{ed}

\title{Reeb graph and quasi-states on the two-dimensional torus}
\author{Frol Zapolsky}

\begin{document}

\maketitle

\begin{abstract}This note deals with quasi-states on the two-dimensional torus. Quasi-states are certain quasi-linear functionals (introduced by Aarnes) on the space of continuous functions. Grubb constructed a quasi-state on the torus, which is invariant under the group of area-preserving diffemorphisms, and which moreover vanishes on functions having support in an open disk. Knudsen asserted the uniqueness of such a quasi-state; for the sake of completeness, we provide a proof. We calculate the value of Grubb's quasi-state on Morse functions with distinct critical values via their Reeb graphs. The resulting formula coincides with the one obtained by Py in his work on quasi-morphisms on the group of area-preserving diffeomorphisms of the torus.
\end{abstract}


\renewcommand{\labelenumi}{(\roman{enumi})}

\section{Introduction and formulation of the result}\label{section_intro}

Following Aarnes \cite{Aarnes_quasi-states}, we give the following definition:
\begin{defin}
If $Z$ is a compact (Hausdorff) space, let $C(Z)$ denote the Banach algebra of all real-valued continuous functions on $Z$. Denote by $C(F)$ the closed subalgebra of $C(Z)$ generated by $F$, that is $C(F)=\{\phi\circ F\,|\,\phi\in C(\im F)\}$. A functional $\eta \fc C(Z) \to \R$ is called a quasi-state if it satisfies
\begin{enumerate}
\item $\eta(1) = 1$;
\item $\eta(F) \geq 0$ for $F \geq 0$;
\item for each $F \in C(Z)$ the restriction $\eta|_{C(F)}$ is linear.
\end{enumerate}
\end{defin}

\noindent Note for further use the nontrivial fact that a quasi-state is monotone, that is $\eta(F)\leq \eta(G)$ for $F\leq G$, and consequently it is Lipschitz with respect to the uniform norm. In particular, it is continuous in the $C^0$-topology on $C(Z)$.

Consider the two-dimensional torus $X=\T^2$, endowed with an area form $\omega$ of total area $1$. In \cite{Grubb_irr_part_constr_tms} Grubb constructed a quasi-state $\zeta$ on $X$, which is invariant under all diffeomorphisms of $X$ preserving $\omega$, and which moreover vanishes on all functions with support contained in an open disk. In \cite{Knudsen_qm_torus} Knudsen asserted the uniqueness of such a quasi-state; for the sake of completeness, we provide a proof below. Our purpose in this note is to calculate the value of Grubb's quasi-state $\zeta$ on a Morse function $H\in C^\infty(X)$ with distinct critical values in terms of its Reeb graph.

\begin{defin}The Reeb graph $\Gamma_H$ of $H$ is obtained from $X$ by collapsing connected components of level sets of $H$ into points.
\end{defin}

\noindent We let $\pi \fc X \to \Gamma_H$ denote the quotient map. It can be easily seen that the graph $\Gamma_H$ has a unique circular subgraph $\Gamma$ with vertices $s_1,\dots,s_k$, and pairwise disjoint trees $T_1,\dots,T_k \subset \Gamma_H$, such that each $T_j$ intersects $\Gamma$ precisely at $s_j$, and $\Gamma_H = \Gamma \cup \bigcup_j T_j$. For $j=1,\dots,k$ let $\alpha_j := H(s_j)$. Then
\begin{thm}\label{thm_main}
$$\zeta(H) = \int_XH\omega+\sum_{j=1}^k \int_{\pi^{-1}(T_j)}(\alpha_j-H)\omega\,.$$
\end{thm}

\noindent The proof is given below. What is surprising is the fact that exactly the same formula (up to a factor of $2$) was obtained by Py in \cite{Py_calabi_qm_on_torus} for the value on $H$ of the functional $\theta$ stemming from a quasi-morphism that he defined on the identity component of the group of symplectomorphisms of $(X,\omega)$, provided that $H$ satisfies the additional requirement of having zero mean, that is $\int_XH\omega = 0$. Denote by $C^\infty_0(X)$ the set of smooth functions on $X$ having zero mean. Then we obtain
\begin{coroll}\label{coroll_Py_eq_Grubb}Py's functional is the restriction to $C^\infty_0(X)$ of a quasi-state, namely, Grubb's quasi-state.
\end{coroll}

\noindent The corollary follows from the discussion, except one technical point, which is accounted for below.

\begin{rem}This note arose as an attempt to understand Py's functional on the torus and to extend the results of Rosenberg \cite{Rosenberg_quasi_states_calabi_qm_h_genus} on surfaces of higher genus to the case of the torus.
\end{rem}

\begin{acknow}I would like the thank Prof.\ Leonid Polterovich who suggested that I publish this result, and provided valuable remarks.
\end{acknow}

\section{Proofs}

In order to carry out the calculation of theorem \ref{thm_main}, we need another definition from \cite{Aarnes_quasi-states}.
\begin{defin}For a compact space $Z$ let $\cA$ denote the collection of subsets of $Z$ which are open or closed.
A function $\tau\fc\cA\to[0,1]$ is called a topological measure if it satisfies
\begin{enumerate}
\item $\tau(Z) = 1$;
\item $\tau(Z-K)+\tau(K)=1$ for compact $K \subset Z$;
\item $\tau(K\cup K')=\tau(K)+\tau(K')$ for disjoint compact $K,K' \subset Z$;
\item $\tau(K)\leq \tau(K')$ for compact $K,K'\in\cA$ such that $K\subset K'$;
\item for open $U\subset Z$, $\tau(U)=\sup\{\tau(K)\,|\,K\subset U,\,K\text{ compact}\}$.
\end{enumerate}
\end{defin}
Quasi-states are in bijection with topological measures -- this is Aarnes's representation theorem \cite{Aarnes_quasi-states}. In one direction, if $\tau$ is a topological measure, then the value of the corresponding quasi-state $\eta_\tau$ on a function $F\in C(X)$ is
$$\eta_\tau(F)=\max_XF - \int_{\min_XF}^{\max_XF}b_F(t)\,dt\,,$$
where $b_F(t) = \tau(\{F\leq t\})$.

Grubb actually constructed a topological measure and the aforementioned quasi-state is, of course, provided by Aarnes's representation theorem. In what follows Grubb's topological measure is denoted by $\tau$. Maintain the notations of section \ref{section_intro}. Since the function in question, $H$, has only finitely many critical values, it suffices to compute $b_H(t) = \tau(\{H\leq t\})$ for a regular value $t$ of $H$. Let $K_j := \pi^{-1}(s_j)$,  $D_j:=\Int\pi^{-1}(T_j)$ for $j=1,\dots,k$ and $S:=X - \bigcup_j D_j = \pi^{-1}(\Gamma)$. Note that the $K_j$ are figures-eight and the $D_j$ are open disks. Denote by $|\cdot|$ the Lebesgue measure corresponding to $\omega$. For $t \in \R$ let $X^t := \{H\leq t\}$.
\begin{lemma}For a regular value $t$ of $H$ we have
\begin{equation}\label{eqn_tm_of_sublevel_set}b_H(t) = |X^t \cap S|+\sum_{j:\alpha_j>t}|D_j|\,.\end{equation}
\end{lemma}

Assuming the lemma for the moment, we proceed to the
\begin{prf}[of theorem \ref{thm_main}]
Without loss of generality, $\alpha_1 < \dots < \alpha_k$. Denote also $M = \max_X H$ and $m=\min_X H$. Then we have
\begin{align*}\zeta(H)&=M-\int_m^M b_H(t)\,dt = M - \int_m^M\Big[|X^t\cap S| + \sum_{j:t>\alpha_j}|D_j|\Big]\,dt\\
&=M - (M-\alpha_k)|S|- \int_{\alpha_1}^{\alpha_k}|X^t\cap S|\,dt\\
&\qquad\qquad-\sum_{j=1}^{k-1}(\alpha_{j+1}-\alpha_j)\Big|\bigcup_{i=1}^{j} D_i\Big|-(M-\alpha_k)\sum_{j=1}^k|D_j|\\
&=M - M|S|+\alpha_k|S|-M(1-|S|) - \int_{\alpha_1}^{\alpha_k}|X^t\cap S|\,dt + \sum_{j=1}^k\alpha_j|D_j|\\
&=\alpha_k|S| - \int_{\alpha_1}^{\alpha_k}|X^t\cap S|\,dt + \sum_{j=1}^k\alpha_j|D_j|=\int_S H\omega + \sum_{j=1}^k\alpha_j|D_j|\\
&=\int_X H \omega + \sum_{j=1}^k\Big(\alpha_j|D_j| - \int_{D_j}H\omega\Big) = \int_X H \omega + \sum_{j=1}^k \int_{D_j}(\alpha_j-H)\omega\,,
\end{align*}
as required. Here we used the fact that if $(Y,\sigma)$ is a finite measure space, $G$ is a bounded measurable function on $Y$ with maximum $\beta$ and minimum $\alpha$, then
$$\int_Y G\,d\sigma = \beta\sigma(Y) - \int_\alpha^\beta \sigma(\{G\leq t\})\,dt\,;$$
here we have $Y=S$, $\sigma=|\cdot|$, $G=H|_S$, $\alpha = \alpha_1$, $\beta=\alpha_k$. \qed
\end{prf}

\begin{prf}[of the lemma]
In order to establish formula \eqref{eqn_tm_of_sublevel_set}, we need to compute $\tau(\{H\leq t\})$. The topological measure $\tau$ can be described as follows. A topological measure is completely determined by its values on compact submanifolds (with boundary) of full dimension, in our case, on compact subsurfaces, see \cite{Zapolsky_isotopy_invariant_tms_surf_high_genus}. If $W \subset X$ is a compact subsurface, let $\widehat W \subset X$ denote the unique subsurface with no contractible boundary components such that (i) $\partial \widehat W \subset \partial W$, (ii) $\ol{\Int W \cap \Int \widehat W} \supset \partial \widehat W$, and (iii) $\ol{W \triangle \widehat W}$ is contained in a (finite) disjoint union of closed disks. Then $\tau(W) = |\widehat W|$. The uniqueness of $\widehat W$, as well as certain properties of the map $W\mapsto \widehat W$ are established in \cite{Zapolsky_isotopy_invariant_tms_surf_high_genus}; the notation $\ed W$ is used there instead of $\widehat W$.

Therefore we need to find $\widehat {X^t}$ for a regular value $t$ of $H$. If $\gamma \subset X$ is a simple closed contractible curve, then it bounds a unique closed disk, which we denote $D(\gamma)$.
\begin{defin}Let $W \subset X$ be a subsurface. A contractible boundary component $\gamma \subset \partial W$ is called exterior (with respect to $W$) if $\ol{\Int D(\gamma) \cap \Int W} \supset \gamma$ and interior otherwise. A contractible boundary component $\delta \subset \partial W$ is called maximal if for any other contractible boundary component $\delta'$ we have either $D(\delta')\subset D(\delta)$ or $D(\delta')\cap D(\delta) = \varnothing$.
\end{defin}
The following lemma is proved in \cite{Zapolsky_isotopy_invariant_tms_surf_high_genus}:
\begin{lemma}\label{lemma_eliminating_disks}Let $W\subset X$ be a subsurface, and assume that $\gamma_1,\dots,\gamma_r$ are the maximal exterior boundary components of $W$, while $\delta_1,\dots,\delta_s$ are the maximal interior boundary components of $W$. Then
\[\widehat W = \Big(W \cup \bigcup_i D(\delta_i) \Big) - \bigcup_j D(\gamma_j)\,.\]\qed
\end{lemma}

Thus in order to compute $\widehat {X^t}$, we need to know the maximal boundary components of $\partial X^t = H^{-1}(t)$, and also which components are exterior and interior with respect to $X^t$. Note that if $\gamma \subset H^{-1}(t)$ is a contractible curve, then it must be contained in one of the disks $D_j$. Before determining the exterior and interior maximal components of $X^t$, we need an auxiliary result, which also appears in \cite{Zapolsky_isotopy_invariant_tms_surf_high_genus}:
\begin{lemma}Let $W \subset X$ be a subsurface contained in an open disk. Then any boundary component of $W$ is contractible and any maximal boundary component of $W$ is exterior with respect to $W$. \qed
\end{lemma}

\begin{lemma}\label{lemma_exterior_interior_X_t}Let $\gamma \subset H^{-1}(t)=\partial X^t$ be a maximal boundary component. Then it is exterior with respect to $X^t$ if and only if $t < \alpha_j$ where $j$ is the unique index such that $\gamma\subset D_j$.
\end{lemma}
\begin{prf}Let $W = D_j \cap X^t$. If $\alpha_j > t$, then $X^t$, and hence also $W$, is disjoint from $K_j$, which contains the topological boundary of $D_j$. Hence $W$ is a compact subsurface of $X$ contained in an open disk. Since $\gamma$ is a maximal boundary component of $X^t$, it is a maximal boundary component of $W$ and hence is exterior with respect to $W$, by the previous lemma, therefore it is exterior with respect to $X^t$.

To obtain the other direction, repeat the argument for $\{H\geq t\}$ and note that a maximal component $\delta \subset H^{-1}(t) = \partial X^t = \partial \{H \geq t\}$ is exterior with respect to $X^t = \{H \leq t\}$ if and only if it is interior with respect to $\{H \geq t\}$. \qed
\end{prf}

It follows from lemmas \ref{lemma_eliminating_disks}, \ref{lemma_exterior_interior_X_t} that if $\alpha_j > t$, then the disk $D_j$ is disjoint from $\widehat{X^t}$ while if $\alpha_j < t$, then the disk $D_j$ is entirely contained in $\widehat {X^t}$. Consequently
\[\widehat{X^t} = (X^t\cap S)\cup\bigcup_{j:t>\alpha_j}D_j\,,\]
and so
\[b_H(t) = \tau(X^t) = |\widehat{X^t}| = |X^t \cap S| + \sum_{j:t>\alpha_j}|D_j|\,,\]
as asserted. \qed
\end{prf}

\begin{prf}[of corollary \ref{coroll_Py_eq_Grubb}]
So far we have shown that Py's functional $\theta$ and Grubb's quasi-state $\zeta$ coincide on the subset of $C^\infty_0(X)$ consisting of Morse functions with distinct critical values. Since $\zeta$ is a quasi-state, it is continuous in the $C^0$-topology. Using the techniques of section 8 of \cite{Rosenberg_quasi_states_calabi_qm_h_genus}, it is possible to show that $\theta$ is continuous in the $C^2$-topology. The set of Morse functions with distinct critical values is $C^2$-dense in $C^\infty_0(X)$, \emph{a fortiori} $C^0$-dense, hence we are done. \qed
\end{prf}

We now turn to the aforementioned uniqueness. Grubb's topological measure is easily seen to be invariant under any symplectomorphism. We now show
\begin{prop}Let $\sigma$ be a topological measure on $X$ which vanishes on disks and is invariant under symplectic isotopies. Then $\sigma$ equals Grubb's topological measure $\tau$.
\end{prop}

\begin{prf}Since $\sigma$ vanishes on disks, $\sigma(W) = \sigma(\widehat W)$ for any subsurface $W \subset X$. Thus it suffices to show that if $W$ is a subsurface with no contractible boundary components, then $\sigma(W) = |W|$. Using additivity, it is enough to show this in case $W$ is a closed smoothly embedded annulus with non-contractible boundary circles.

Choose a system of coordinates $(p,q)$ on $X$ such that $\omega = dp \wedge dq$. Call an embedded non-contractible circle linear if it has a parametrization of the form $\R/\Z \ni t \mapsto (p_0+kt,q_0+lt)$, where $k,l$ is a pair of mutually prime integers, the slope of the circle. A closed smoothly embedded annulus is linear if its boundary circles are (in this case they must have the same slope). By a standard, thought somewhat lengthy, argument, it can be shown that any closed smoothly embedded annulus can be symplectically isotoped to a linear one of the same area. We shall show that the value of $\sigma$ on any such annulus equals its area.

Clearly translations are symplectic isotopies so that all the linear annuli of the same slope and the same area are symplectically isotopic and so share the same value of $\sigma$.

Fix a slope; in what follows all annuli will be linear and have the chosen slope. Consider an annulus $W$ of area $\frac 1 n$ where $n \geq 2$ is a natural number. Since we can fit $n-1$ pairwise disjoint translates of $W$ into $X$, we have, by the additivity and monotonicity of $\sigma$, $(n-1)\sigma(W) \leq 1$ so that
$$\sigma(W) \leq \frac 1 {n-1}\,.$$
Now consider an annulus $W$ of area $\frac 1 {2N}$ where $N$ is a natural number and let $O = \Int W$. The whole torus $X$ is the disjoint union of $N$ translates of $W$ and $N$ translates of $O$, so we obtain $N(\sigma(W)+\sigma(O)) = 1$, whence
$$\sigma(W) \geq \sigma(O) = \frac 1 N - \sigma(W) \geq \frac 1 N -\frac 1 {2N-1}\,,$$
by the previous inequality.

Finally, let $A$ be an annulus of arbitrary area $\alpha \in (0,1)$. For any natural $N$ let $m_N = \lfloor 2N\alpha \rfloor$. Note that $\lim_{N\to\infty}\frac {m_N}{2N} = \alpha$. We can fit at least $m_N-1$ translates of an annulus of area $\frac 1 {2N}$ into $A$. This means that
$$\sigma(A) \geq (m_N-1)\left(\frac 1 N - \frac 1 {2N-1}\right)\,.$$
On the other hand, we can fit at least $2N - 1 - m_N$ translates of an annulus of area $\frac 1 {2N}$ into the complement $X-A$, so that
$$\sigma(A) = 1 - \sigma(X-A) \leq 1 - (2N - 1 - m_N)\left(\frac 1 N - \frac 1 {2N-1}\right)\,.$$
It follows easily from the last two inequalities that for large $N$
$$\left|\sigma(A) - \frac {m_N}{2N}\right| \leq \frac{\const} N\,.$$
Therefore $\sigma(A) = \lim_{N\to\infty}\frac{m_N}{2N} = \alpha$, as desired. \qed

\end{prf}

\end{document}